\documentclass[11pt]{article}
  \usepackage[all]{xy}
   \usepackage[ansinew]{inputenc}
    \usepackage[psamsfonts]{amssymb}
  \newcommand{\Z}{\mathrm{Z}}

   \def\build#1_#2^#3{\mathrel{\mathop{\kern0pt#1}\limits_{#2}^{#3}}}

    \def\bk{{l\!k}} 

    \def\hH{{H\!H}} 

    \def\hC{{\bf C}} 

   \def\cT{{T\! C}} 

    \def\bB{I\!\!B} 

    \def\bH{I\!\!H} 

  \newcommand{\bBnu}{\widetilde{\bB}} 

 \newcommand{\Omeganu}{\widetilde{\Omega}}

    \title{Duality in Gerstenhaber algebras}

    \author{Yves F\'elix, Luc Menichi and  Jean-Claude Thomas }

\begin{document}

    \maketitle

    \begin{abstract} Let $C$ be a differential graded coalgebra, $
      \overline\Omega C$ the Adams cobar construction and $C^\vee$
      the dual algebra. We prove that for a large class of coalgebras
      $C$ there is a natural isomorphism of Gerstenhaber algebras
      between the Hochschild cohomologies  $HH^\ast   (C^\vee   ,C
      ^\vee )$ and $H\!H^\ast   ( \overline\Omega C ; \overline\Omega
      C )$. This result permits to describe a Hodge decomposition of
      the loop space homology of a closed oriented manifold,  in the
      sense of Chas-Sullivan,  when the
      field  of coefficients is of characteristic zero. 
  \end{abstract}

    \vspace{5mm}\noindent {\bf AMS Classification} : 
16E40, 17A30, 17B55, 81T30, 55P35.

    \vspace{2mm}\noindent {\bf Key words} :  Hochschild cohomology,
    Gerstenhaber algebra, free loop space, loop homology.

  \section{Introduction }

    In the last two decades there  has been  a great deal of interest in
    Gerstenhaber algebras since they arise, as BV-algebras, in
    the BRST theory of topological field    theory and in operad setting as well
    as in string theory.  Let us denote by ${\bf C}^\ast (A;A)$ (resp.
    $H\!H^*(A;A)$)  the Hochschild complex (resp. cohomology) of the
 differential graded
  algebra $A$ with coefficients in  $A$. It is  well
    known that $H\!H^*(A;A)$ is a Gerstenhaber algebra.
A new
   geometrical impact to   Hochschild cohomology  has been given recently by the 
work of Chas and Sullivan
\cite{CS}.

Let $\bk$ be an integral domain 
and let  $V^{\vee}$  denotes the graded dual  of the
graded module  $V$. For a supplemented graded
coalgebra
 $C =\bk \oplus  \overline C $,
$C^\vee$ is a  supplemented   graded algebra (without any
 finiteness  restriction),
and the reduced coproduct $\overline \Delta : \overline C \to \overline 
C \otimes \overline C $ can be iterated
unambiguously to produce $\overline \Delta ^{(k)}$.  If for each
 $x\in \overline C
$ and some $k $, $\overline \Delta ^{(k)} x = 0 $,  the coalgebra $C$
is called  {\it locally  
conilpotent}.  

Our first  result reads

    \vspace{5mm}
    \noindent{\bf Theorem 1.} {\it Let (C,d) be a $\bk$-free
locally  conilpotent  differential graded
    coalgebra and  $\overline\Omega
     C$ the normalized cobar construction on  $C$.
There exists  a homomorphism of 
differential graded  algebras
  $${\cal D} : {\bf C}^\ast (\overline\Omega C;\overline\Omega C)
  \longrightarrow  {\bf 
C}^\ast ( C^\vee
 ;C^\vee )$$
which induces a homomorphism  of Gerstenhaber  algebras 
$$H({\cal D}) : H\!H^\ast (\overline\Omega C;\overline\Omega C) \longrightarrow
H\!H^\ast ( C^\vee
 ;C^\vee ).$$

a) The homomorphism $H(\mathcal{D})$ is an isomorphism
whenever  $C$  is finitely generated in each degree, and  either
 i) $ C =
C^{\geq 0}$ or else\quad 
ii) $   \overline{C}=\overline{C}_{\geq   2}$.

b) If $\bk$ is a field then
$H({\cal D})$ is natural with respect to
  quasi-isomorphisms of differential graded coalgebras
and is an isomorphism  
whenever  $ H(C)$  is finitely generated in each degree, and  either
i) $ C =
C^{\geq 0}$ and $ H^0(C,d) =
\bk$, or else\quad 
ii) $   \overline{C}=\overline{C}_{\geq   2}$.}

   \vspace{3mm}

The proof of Theorem 1   relies  
heavily on  
properties of   free models and on  
 differential graded Lie algebras of
derivations. In order to minimize the periquisites, we do 
not introduce the twisting cochain
formalism. Nonetheless  it   underlies most of the  computational peculiarities 
we perform by hands.

Let us precise  that a {\it free model}  is a  tensor algebra $TV$
together with a
differential $d$ such that  $V$ is the union 
$\displaystyle V = \bigcup _{k=0} ^\infty V(k) $ of an increassing familly  
$V(0)\subset V(1)\subset ... $ such that $d(V(0))=0$ and $ d(V(k)) \subset 
T(V(k-1))$. Observe 
(\cite{FHT}-Proposition 3.1)
that for any differential graded algebra $A$ there is a
quasi-isomorphism of differential graded algebras  
 $(T(V), d) 
\stackrel{\simeq }{\twoheadrightarrow }A$.  If  $C$ is a locally
conilpotent  coalgebra, the normalized   cobar construction on $C$  is
a  free model (5.2).

 We denote by ${\rm Der}A$ the differential graded Lie algebra of derivations on 
$A$ with the commutator 
 bracket $[-,-]$ and 
 differential 
 $D= [d,-]$. We introduce the following extension, denoted  $\widetilde {\rm Der} A$,  of this  differential graded 
Lie algebra:
 $$ 
 \widetilde {\rm Der} A = {\rm Der} A \oplus sA \,, \quad 
 \left\{ 
 \begin{array}{ll} \widetilde D ( \theta + sx) &= D(\theta) - {\rm ad}_x - sd(x) 
\\

{\rm ad}_x(y)&= xy - (-1)^{|x||y|}
yx\\

[ \theta, \theta' +sx] & = [\theta, \theta'] + (-1)^{|\theta|} s\theta( x)\\

[sx,sy] &= 0 
 \end{array} 
 \right. 
 \,\mbox{with}\left\{ 
 \begin{array}{ll} 
 \theta \,, \theta ' \in {\rm Der} A\\ 
 x\,, y \in A 
 \end{array} 
 \right. 
 $$
where $(sA)_i =A_{i-1}$. The next result  is a corner stone in the proof of Theorem 1.

 \vspace{3mm}\noindent {\bf Theorem 2.} 
 {\it Let $A=(TV,d)$ be a free model. There exist  two  injective quasi-isomorphisms  of differential graded
Lie algebras
$$
s{\hC}^\ast (A;A) \stackrel{\simeq}{\longleftarrow} \widetilde {\rm Der} A 
\stackrel{\simeq}{\longrightarrow}  {\rm Der} \widetilde A
$$
where $ \widetilde A$ is the differential graded algebra   $ \widetilde A 
=(T(V\oplus \bk \epsilon ), \tilde d)$ with $\tilde d
\epsilon =
\epsilon ^2
$ and
$
\tilde d v = dv +
\epsilon v - (-1)^{|v|} v \epsilon $, $v \in V$.}
\vspace{3mm}

If we introduce, following  Getzler and Jones \cite{G:Cartan_Homotopy,GJ}, the
  Hochschild complex, $\hC _\infty^\ast (E;E)$ of an $A_\infty$-coalgebra
$E$,  Theorem 2 can be
read: {\it There exist  injective quasi-isomorphisms of differential
  graded Lie algebras
$$
s{\hC}^\ast (A;A) \stackrel{\simeq}{\longleftarrow} \widetilde {\rm Der} A 
\stackrel{\simeq}{\longrightarrow} 
 \hC _\infty^\ast (\bk\oplus sV;\bk\oplus sV).
$$}

Since the graded Lie algebra $s{\hH}^\ast (A;A)$ is
natural with respect to quasi-isomorphisms (3.4) we
deduce from Theorems 1 and 2:

 \vspace{3mm}
    \noindent{\bf Corollary 1.} { \it Consider $C$ as in Theorem
      1. Let
 $\varphi : (TV,
  d)\stackrel{\simeq} \to
   C^{\vee} $
   and $\psi : (TW , d)\stackrel{\simeq} \to \overline{\Omega} C $
    be free  models, then  there exists  isomorphisms of
      graded Lie algebras:
    $$ H\left( \widetilde{\rm Der} (TW,d )\right)\,
\cong\,
 sH\!H^\ast  (C^{\vee};C^{\vee})
\cong
H\left(\widetilde {\rm Der } (TV,d)\right).
$$
}

\vspace{1mm}
The geometric meaning of Theorem 1 is in terms of loop homology.
  Recall that the loop homology of a closed
   orientable
    manifold $M$ of dimension $d$  is the ordinary homology
    of the free loop space $M^{S^1}$ with degrees shifted by $d$, i.e.
  $\mathbb
    H_\ast  (M^{S^1}) = H_{\ast  +d}(M^{S^1})$. In \cite{CS}, Chas and
  Sullivan
     have defined a   Batalin-Vilkovisky 
algebra structure on $\mathbb
    H_\ast  (M^{S^1})$  . Thus    a loop product and a loop
 bracket are defined  such that   $\mathbb H (M^{S^1}) $ is a
Gerstenhaber algebra.  More recently, in his thesis \cite{Tr},
Thomas Traddler has introduced the notion  
of $A_\infty$-bimodules and $\infty$-inner product
for $A_\infty$- algebras such that 

\noindent  a) the Hochschild homology   is defined and is naturally a BV-algebra 
for an 
 $A_\infty$- algebra with $\infty$-inner product 

\noindent  b)  the geometric chains on a manifold naturally possesses via 
intersection 
 the structure of an $ A_\infty$-algebras with  $\infty$-inner product 

\noindent  c) for simply connected closed  manifolds a)  and b) are compatible with the 
"Chas-Sullivan" BV-structure.

     Denote by    ${\cal C}_\ast
  (M)$
   (resp.
   ${\cal C}^\ast  (M)$) the
    normalized singular chain coalgebra (resp. cochain algebra) on $M$ with
    coefficients in  $\bk$.
    The cap product with a fixed fundamental cycle, ${\cal C}^\ast  (M) \to
  {\cal
    C}_{d-\ast  }(M)$, is  not a map
    of ${\cal C}^\ast  (M)$-bimodules. Nonetheless, (\cite{FTV} Theorem 2),
   Poincar\'e
    duality induces an isomorphism $
     {\cal P} : H\!H^\ast  ({\cal C}^\ast  (M); {\cal C}_\ast  (M))
   \stackrel{\cong}{\longrightarrow}
    H\!H^{\ast  -d}({\cal C}^\ast  (M); {\cal C}^\ast  (M))\,.
    $
     The Jones isomorphism \cite{J}
    $
    J :  H_{n}(M^{S^1} )  \stackrel{\cong}{\longrightarrow}
    H\!H^n ({\cal C}^\ast (M); {\cal C}_\ast (M))
    $
    composed with ${\cal P}$, is an isomorphism of graded modules. 
  R. Cohen and J. Jones (\cite{CJ}-corollary 10) have announced  that this isomorphism also identify   the  loop product on $\mathbb
    H_\ast  (M^{S^1})$    with
    the Gerstenhaber product on ${H\!H}^{\ast   } ( C^\ast (M); C^\ast (M) )$.
Then Theorem 1 and naturality of ${H\!H}^{\ast   }(A;A)$ (3.4) would imply a ``Gerstenhaber algebra analog'' of a
  result obtained by  Burghelea-Fiederowicz \cite{BF} and
    Goodwillie \cite{Go}: {\it   Let  $M$ be   a simply-connected closed
    oriented   manifold. The loop algebra
    $\bH_\ast  (M^{S^1})$  is isomorphic to the Gerstenhaber 
    algebra
    $H\!H^\ast({\cal C}_\ast  (\Omega M), {\cal C}_\ast  (\Omega M)$.}
  Here $\Omega X$ denotes the based loops on a pointed space
  $X$. Indeed for any pointed 1-connected space $X$ there is a natural
  equivalence $\overline \Omega \mathcal C^\ast(X) \to \mathcal
  C^\ast(\Omega X)$, \cite{ACC}.

 \vspace{3mm}

Corollary 1 is particulary interesting when   $\psi : (T(W),d)\to {\cal 
C}_*(\Omega M)$
is an
  Adams-Hilton
    model \cite{AH}  for $M$, since the Adams-Hilton model of a space is
completely determined by a
    cellular decomposition of $M$. If  $\bk$ is a field of
    characteristic zero we can go further.   The
Adams-Hilton
 model of $M$ is   the universal  enveloping algebra of a differential
 graded Lie algebra $(L,d) = (\mathbb L(V),d)$ \cite{An}. Then,
$$
\begin{array}{lll}
\mathbb H_{*}(M^{S^1})
&\cong \hH^*({\cal C}_*(\Omega M);{\cal C}_*(\Omega M)) &\mbox{(Corollary
  1)}\\
&\cong\hH^*(U(L,d);U(L,d)) & \mbox{(by naturality) }\\
&\cong  \mbox{Ext}_{{U(L,d)}^e}( U(L,d),U(L,d)) &  \\
&\cong\mbox{Ext}_{{U(L,d)}}( \bk ,U(L,d)) & \mbox{(\cite{CE},
 Theorem X$\!$I$\!$I$\!$I.5.1), }
\end{array}
$$
where   $UL$ is considered as an $UL$-module via the adjoint
 representation. Denote by $\Gamma^n(V)$ the vector space generated by
 the elements $\sum_{\sigma\in \Sigma_n} v_{\sigma (1)}\cdots
 v_{\sigma (n)}$, with $v_i \in V$. Then the vector spaces $\Gamma^n(V)$
 are stable under  the adjoint action of $UL$ and $UL \cong
 \oplus_{n\geq 0} \Gamma^n(V)$.
This gives the  ``Hodge decomposition'':

\vspace{3mm}
\noindent{\bf Corollary 2.} {\it Under the above hypothesis there
exists
isomorphisms of graded vector spaces:
 $$
\mathbb H_{*}(M^{S^1}) =\mbox{Ext}_{{U(L,d)}}( \bk ,U(L,d))= \oplus_{n\geq 0}
{\rm Ext}_{UL}(\bk,
 \Gamma^n(V))\,.$$
}

   \vspace{5mm}
   The remaining of the  paper is organized as follows:

  \noindent  2) Sketch of the proof of Theorem 1.

  \noindent  3)  Hochschild cohomology of a differential graded algebra.

\noindent 4) Proof of
   Theorem 2.

\noindent 5) The Hochschild cochain complex of a differential
graded
 coalgebra.

   \noindent 6) Proof of Propositions  B and C (see below).

\noindent\section{ Sketch of the Proof of Theorem 1.}

\vspace{5mm}
 \noindent {\bf 2.1. Notation.}  
In the rest of the paper, except in 2.6 and 2.7, $\bk$ will be a
principal ideal domain.
     If $V =\{V_i\}_{i \in \mathbb Z}  $ is a
(lower)
 graded $\bk$-module  (when we need upper graded $\bk$-module  we put
 $V_i= V^{-i}$ as usual) then:

    \vspace{1mm}\noindent a) $(sV)_n =V_{n-1}\,, \,\, (sV)^n =V^{n+1}$,

    \vspace{1mm}\noindent b)  $TV $ denotes the tensor algebra on $V$,
while we denote by $T\!C(V)$ the free
 supplemented coalgebra generated by $V$.

 Since we work with
    graded differential objets, we will make a
    special attention to signs.  Recall that if
  $M=\{M_i\}_{i \in \Z}$ and  $N= \{N_i\}_{i \in \Z}$ are differential
graded $\bk$-modules  then

 \noindent  a) $M\otimes N$ is a differential graded $\bk$-module : $(M\otimes N)_r =
\oplus _{p+q=r} M_p
  \otimes N_q$, $ d_{M\otimes N}= d_M \otimes id_N + id_M
  \otimes d_N$,

 \noindent  b) $Hom(M,N)$ is a differential graded $\bk$-module  : $Hom_n(M,N)=
\prod_{k-l=n} Hom(M_l,N_k)$,
  $Df  = d_M\circ f - (-1) ^{|f|} f\circ d_N$,

\noindent c) the commutator, $[f,g]=f\circ g-(-1)^{\vert f\vert\vert
  g\vert}g\circ f$, gives to the differential
graded $\bk$-module
$\mbox{End}(M)= \mbox{Hom}(M,M)$ a structure of differential graded Lie algebra,

\noindent d) if $C$ is a differential graded coalgebra with diagonal
$\Delta$
 and $A$ is a differential graded algebra with product $\mu$ then the
 cup product,
$f\cup g= \mu\circ (f\otimes g) \circ \Delta$, gives  to the differential
graded $\bk$-module
$\mbox{Hom}(C,A)$ a structure of differential graded  algebra.

\vspace{5mm}
 \noindent {\bf 2.2.}  A {\it (graded) Gerstenhaber algebra} is a
 commutative  graded algebra $G=\{G\}_{i\in
\Z}$ with a degree 1  linear map 
$$
G_i \otimes G_j \to G_{i+j+1} \,, \quad x\otimes y \mapsto \{x,y\}
$$ 
such that:

\noindent a) the suspension of $G$  is a   graded Lie
 algebra  with bracket
$$
(sG)_i \otimes (sG)_j \to (sG)_{i+j} \,,\quad sx\otimes sy \mapsto [sx,sy]:=s 
\{x,y\}
$$

\noindent b)  the product  is   compatible with the bracket,  $\{-,-\}$. 

This last condition means that for any  $a  \in G_k$  
the adjunction map $ad_a :G_i \to G_{i+k+1} \,, \quad
b\mapsto \{a,b\}$ is  a
$(k+1)$-derivation: ie. for $a$, $b$, $c\in G$,
$\{ a, b c\} = \{ a,b\} c + (-1)^{\vert b\vert (\vert a\vert +1)} 
b \{ a,c\}$.

 A {\it homomorphism of  Gerstenhaber algebras} $f:G\to G'$
 is a 
homomorphism of
graded algebras such that $sf : sG \to sG'$ is a homomorphism of graded Lie 
algebras.

For our purpose it is convenient to introduce the notion of {\it
  pre-Gerstenhaber algebra}. This a differential graded algebra
  $\mathcal G=(\{\mathcal G \}_{i}, d) $ together with a degree 1 homomorphism of differential
  graded modules $ \alpha :{ \mathcal G}\to {\mathcal L}$ where
  ${\mathcal L}=(\{\mathcal L \}_{i}, D )$ denotes a differential graded Lie algebra. Observe that

\noindent a) $D\circ \alpha =-\alpha \circ d $.

\noindent b) $\alpha $ induces on $s{\mathcal G}$ a structure of
  graded a  Lie algebra compatible with the differential $sd$ while no
  compatibility condition with the product is required.

A {\it  homomorphism of pre-Gerstenhaber algebras} is a 
commutative diagram of
homomorphisms of differential  graded modules 
$$
\begin{array}{lclcl}
s{\mathcal G} & \stackrel{\alpha }{\rightarrow }& {\mathcal L}\\
sf \downarrow && \downarrow g\\
s{\mathcal G'}  &  \stackrel{\alpha' }{\rightarrow }& {\mathcal L}\\
\end{array}
$$
such that $f$ is a homomorphism of differential graded algebras and  $g: L \to L'$ 
is a homomorphism of differential graded Lie
algebras. Clearly, $f$ induces a homomorphism of graded algebras $
H(f): H(\mathcal  G) \to H(\mathcal G') $ and a homomorphism of
differential graded Lie algebras  $
sH(f): H(s \mathcal  G) \to H(s \mathcal  G') $. 

If $H(f)$ is an isomorphism of graded modules then $f$ is called a
{\it quasi-isomorphism of pre-Gerstenhaber algebras}. 

If the structure of pre-Gerstenhaber algebra on $\mathcal G$
(resp. $\mathcal G'$ ) induces a structure of Gerstenhaber algebra on
$H(\mathcal G)$ (resp. on $H(\mathcal G')$) then $H(f)$ is a
homomorphism of Gerstenhaber algebras.

\vspace{5mm}
 \noindent {\bf 2.3.}  Let $A$ be a differential graded algebra.
We consider the degree 1  
isomorphism $\beta_A$ that extends a linear map to a coderivation
$$
{\bf C}^\ast (A;A):= \mbox{Hom}( T\!C(sA), A)   
 \stackrel{\beta_A }{\rightarrow }\mbox{Coder } (\widetilde{\mathbb B}(A))
$$
where

\noindent i)   $ \widetilde{\mathbb B} (A)$ is the non-unital  bar construction (3.1),

\noindent ii) $\mbox{Coder } (C) $ is the  differential graded Lie algebra of 
coderivations of  $C$.

The map $\beta_A$ defines the Hochschild complex  ${\bf C}^\ast (A;A)$
as a pre-Gerstenhaber algebra. (It induces on  $\hH^{\ast}(A;A)$ 
the usual Gerstenhaber algebra structure \cite{St}.)  

\vspace{5mm}
 \noindent {\bf 2.4.}  Dually, if $C$ is a supplemented  differential
 graded coalgebra the degree 1  
isomorphism $\gamma_C$ that extends a linear map to a derivation
$$
{\bf C}^\ast (C;C):= \mbox{Hom}(C, Ts^{-1}C)  
 \stackrel{\gamma_C }{\rightarrow } \mbox{Der } (\widetilde \Omega C)
$$
where  $ \widetilde \Omega C$ is the non-counital cobar construction (5.1), makes the 
Hochschild complex of  $C$ into a pre-Gerstenhaber
algebra.

We  also consider the pre-Gerstenhaber algebra $ \overline {\bf C}^\ast 
(C;C) $ together with  a  degree 1 linear isomorphism:
$$
\overline{\bf C}^\ast (C;C):= \mbox{Hom}(C, Ts^{-1}\overline C)  
 \stackrel{\overline{\gamma}_C}{\rightarrow }\widetilde{\mbox{Der }} ( 
{{\overline\Omega}} C)
$$
where 

\noindent i) $ \overline \Omega C$ is the normalized  cobar
construction

\noindent ii) $\widetilde{\mbox{Der }}  (A ) $ is the differential graded Lie algebra  
considered in the introduction.

\vspace{5mm}
 \noindent {\bf 2.5. Intermediate  results.}
 
\vspace{5mm}
\noindent{\bf Proposition A.} {\it Let $C=\bk \oplus \overline C $ be
  a $\bk$-free 
locally conilpotent  differential graded coalgebra. The  inclusion $
\overline C
\hookrightarrow C$  induces, by naturality, a quasi-isomorphism of pre-Gerstenhaber algebras
$$
\overline{\bf C}^\ast (C;C) \to{\bf C}^\ast (C;C)\,.
$$
}

\vspace{3mm}

\noindent{\bf Proposition B.} {\it Let $C $ be a  differential 
graded coalgebra.
The usual linear duality induces an homomorphism of
pre-Gerstenhaber algebras
${\bf C}^\ast (C;C)\longrightarrow {\bf C}^\ast (C^\vee ;C^\vee )$
which is

\noindent i) an isomorphism whenever $C=C^{\geq 0}$
is a free graded $\bk$-module of finite type,

\noindent ii) a quasi-isomorphism if $C=\bk\oplus C_{\geq 2}$
is a free graded $\bk$-module of finite type.}

\vspace{3mm}

\noindent{\bf Proposition C.} {\it Let $C$ be a $\bk$-free locally conilpotent
  differential graded coalgebra.
The bar-cobar adjunction
induces  a quasi-isomorphism of differential graded algebras 
$$  {\bf C}^\ast (\overline 
\Omega C;\overline \Omega  C)
\longrightarrow  
\overline {\bf C}^*(C;C)
$$
 which  admits a linear 
section $ \Gamma$ such that $ s\Gamma : s
\overline {\bf C}^*(C;C) \longrightarrow  s{\bf C}^\ast (\overline
\Omega C;\overline
\Omega  C)$ is a homomorphism of differential graded Lie algebras.}

\vspace{5mm}
 \noindent {\bf 2.6.} Let $C$ be a $\bk$-free locally conilpotent differential
 graded coalgebra. The composite
$$\mathcal D_C : 
{\bf C}^\ast (\overline \Omega C;\overline \Omega  C) 
\stackrel{\mbox{\small Prop. C}}{\longrightarrow}\overline {\bf
C}^*(C;C)\stackrel{\mbox{\small Prop. A}}{\longrightarrow}
 {\bf C}^*(C;C) \stackrel{\mbox{\small Prop. B}}{\longrightarrow}
{\bf C}^\ast (C^\vee ;C^\vee )
$$
is a homomorphism of differential graded algebras but not a
homomorphism of pre-Gerstenhaber algebra. Nonetheless, 
$$
H(\mathcal D_C) :{\hH}^\ast (\overline \Omega C;\overline \Omega  C) 
\longrightarrow
{\hH}^\ast (C^\vee ;C^\vee )
$$
is a homomorphism of Gerstenhaber algebras.
It results from the constructions involved in Propositions A, B
and C that $H(\mathcal D_C)$ is natural with
respect to quasi-isomorphisms of differential graded coalgebras:
Given a quasi-isomorphism of locally conilpotent
supplemented $\bk$-free differential graded coalgebras
$f : C \to D$ such that
either i) $\overline{C}=\overline{C}^{\geq 0}$ and
$\overline{D}=\overline{D}^{\geq 0}$ or else ii) $\overline{C}=\overline{C}_{\geq 2}$ and
$\overline{D}=\overline{D}_{\geq 2}$,
then by Remark 2.3 of~\cite{ACC},
$\overline\Omega(f):\overline{\Omega}C\buildrel{\simeq}\over\rightarrow\overline{\Omega} D$
is a quasi-isomorphism of $\bk$-free differential
graded algebras. Moreover, if $\bk$ is a field,
$f:D^{\vee}\buildrel{\simeq}\over\rightarrow C^{\vee}$ is also a quasi-isomorphism between
$\bk$-free differential graded algebras.
By 3.4, there exist isomorphisms of Gerstenhaber algebras
$\hH(\overline\Omega (f))$ and  $\hH(f^\vee)$ such that the
following diagram commutes

 $$ 
 \begin{array}{ccccc} {\hH}^\ast(\overline\Omega C; \overline\Omega C) & \stackrel{{\hH}^\ast(\mathcal{D_C})}{\longrightarrow }& 
{\hH}^\ast(C^\vee;C^\vee)\\ 

\hspace{-12mm}  {\hH}^\ast(\overline\Omega f) 
\downarrow && 
\hspace{18mm} \uparrow {\hH}^\ast (f^\vee)  \\

 {\hH}^\ast(\overline\Omega D;\overline \Omega D)&  
\stackrel{{\hH}^\ast(\mathcal{D}_D)}{\longrightarrow } &{\hH}^\ast(D^\vee;D^\vee). 
 \end{array}
 $$

\vspace{5mm}
 \noindent {\bf 2.7.}  The   reduction process described below permits
 to deduce part b) of Theorem 1 from part a).

In case ii), let $C$ be a $\bk$-free supplemented differential graded coalgebra
such that $H(C)$ is finitely generated in each degree,
$\overline{C}=\overline{C}_{\geq 2}$ and $H_2(C)$ is $\bk$-free.
By Proposition 4.2 of \cite{ACC}, there exists a free model $TV$
and a quasi-isomorphism of augmented differential graded algebras
$\varphi:TV\buildrel{\simeq}\over\rightarrow C^{\vee}$, where
$V=V^{\geq 2}$ is $\bk$-free of finite type.
We consider then the previous diagram with
$f: C \hookrightarrow C^{\vee\vee}\buildrel{\varphi^{\vee}}\over\rightarrow TV^\vee=D$.  
By  Propositions A, B and C, the  lower line in the above
 diagram is an isomorphism thus so is  the upper line.

The case i) is similar.
\vspace{3mm}

Required definitions and proofs are detailed  in the following sections.

\noindent\section{The Hochschild complex of a differential graded algebra.}

\vspace{5mm}
 \noindent {\bf 3.1.}   Let  $(A, d)$ be a differential graded supplemented 
algebra,
    $ A = \bk \oplus \overline A$, and $(M, d) $ (resp. $(N,d)$) be a  right
 (resp. left) differential graded $ A
    $-bimodule.
    The {\it two-sided bar constructions}, $\bB(M; A; N)$  and
 $\overline{\bB}(M; A; N)$
    are defined as follows:
    $$
    \bB _k (M; A; N) =  M \otimes T^k s  A 
    \otimes N\,, \hspace{1cm}\overline{\bB} _k (M; A; N) =  M \otimes T^k
    s\overline A  \otimes N$$
     A generic element   is written $
    m[a_1|a_2|...|a_k]n$  with  degree   $| m|+ | n|+ \sum_{i=
    1}^k (|s a_i|)\,.$ The differential $d= d_0+d_1$ is defined by
    $$
    d_0  : {\bB} _k  (M;A;N) \to {\bB} _k  (M;A;N)\,, \quad
    d_0  : \overline{\mathrm{B}}_k  (M;A;N) \to \overline{\mathrm{B}}_k  
(M;A;N)\,,$$
    $$  d_1  : \bB _k
    (M;A;N) \to \bB_{k-1}  (M;A;N)\,,\quad d_1  : \overline{\bB}_k
    (M;A;N) \to \overline{\bB}_{k-1}  (M;A;N)\,,
    $$
    with  
    $\epsilon _i = | m| + \sum _{j<i} (|s a_j|)$ and:
    $$
    \renewcommand{\arraystretch}{1.5}
    \begin{array}{rll}
    d_0  ( m[a_1|a_2|...|a_k]n)& = d( m) [a_1|a_2|...|a_k]n  -
   \displaystyle\sum
    _{i=1}^k (-1)^{\epsilon _i}  m[a_1|a_2|...|d(
    a_i)|...|a_k]n\\ &+ (-1) ^{\epsilon _{k+1}}  m[a_1|a_2|...|a_k]d(n)
    \\[2mm]
    d_1(m[a_1|a_2|...|a_k]n)&= (-1) ^{| m|}  ma_1[a_2|...|a_k]n +
    \displaystyle\sum _{i=2}^k (-1) ^{\epsilon _i}
  m[a_1|a_2|...|a_{i-1}a_i|... |
   a_k]n \\
    &- (-1)^{\epsilon _{k}}  m[a_1|a_2|...|a_{k-1}]
    a_k n
    \end{array}
     \renewcommand{\arraystretch}{1}
    $$

 Hereafter we will consider the {\it normalized} and the {\it non-unital}  bar constructions on $A$:
  $$
  \overline{\bB}(A)=\overline{\bB}(\bk;A;\bk)=\left(T\!C(s\overline A), 
\overline d\right)  \,,\quad
 \bBnu(A)=\overline{\bB}(\bk\oplus A) = \left(T\!C(sA) , \tilde d\right) \,.
 $$
  (In the latter  $A$ is considered as a non
  unital algebra (\cite{Qu}-p. 142))). The differentials $\overline d $ and 
$\widetilde d$  are  given by the same formula:
 $$
  \renewcommand{\arraystretch}{1.6}
     \begin{array}{ll}
    d([a_1|a_2|\cdots |a_k]) =& - \displaystyle\sum
    _{i=1}^k (-1)^{\epsilon _i}  [a_1|a_2|...|d(
    a_i)|...|a_k] +
    \displaystyle\sum _{i=2}^k (-1) ^{\epsilon _i}
  [a_1|a_2|...|a_{i-1}a_i|... |
    a_k]\,.
    \end{array}
    \renewcommand{\arraystretch}{1}
 $$
We will also frequently use the twisting cochain of $\bBnu(A)$
$$
\tau_{\bBnu A}:\tilde \mathbb B  A \to A \,, \quad  [a_1\vert ...\vert a_k]\mapsto
\left\{
\begin{array}{ll} 0 &\mbox{ if } k\neq 1\\ a_1 &\mbox{ if } k=
  1\end{array}\right.
$$
and the following result:

 \vspace{5mm}\noindent {\bf Proposition.} \cite[Lemma 4.3]{FHT} {\sl If $A$ 
is a differential graded algebra such that $A$ is a $\bk$-free module  then
the canonical map 
$$
\mathbb B(A,A,A)\to A$$ is a semi-free resolution of $A^e$-modules
(Here $A^{e}:=A\otimes A^{op}$ denotes the envelopping algebra).} 

  \vspace{5mm}
\noindent {\bf 3.2.}  Let $ A$ be a supplemented differential
 graded algebra and $M$ a differential graded $A$-bimodule. The canonical 
isomorphism of graded modules
 $$
\Phi_{A,M}:\mbox{Hom}_{A^{e}}(\bB(A;A;A),M)\rightarrow 
 \mbox{Hom}(\cT(sA),M)\,, \quad  f \mapsto \left( 1_A[a_1|...|a_k]1_A \mapsto 
f([a_1|...|a_k]\right)) 
$$  
carries a differential $D_0+D_1$ on $\mbox{Hom}(\cT(sA),M)$.
More explicitly, if
  $f \in  \mbox{Hom}( T\!C(sA), M) $, we have:
   $$D_0(f)([a_1|a_2|...|a_k])
  = d_{M}(f\left([a_1|a_2|...|a_k])\right) - \sum _{i=1} ^k
  (-1)^{\overline \epsilon_i} f([a_1|...|d_Aa_i|...|a_k])$$  and
  $$\renewcommand{\arraystretch}{1.6}
\begin{array}{ll}
D_1(f)([a_1|a_2|...|a_k])&= 
   - (-1)^{|sa_1|\, |f|}  a_1 f([a_2|...|a_k])
- \sum _{i=2} ^k
  (-1)^{\overline \epsilon_i} f([a_1|...|a_{i-1}a_i|...|a_k])\\
 &\qquad  + (-1) ^{\overline\epsilon _k} f([a_1|a_2|...|a_{k-1}])a_k
  \,,
\end{array}
\renewcommand{\arraystretch}{1}
$$
 where
  $\overline \epsilon _i = |f|+|sa_1|+|sa_2|+...+|sa_{i-1}|$.

 The {\it Hochschild cochain complex} of $A$ with coefficients in the 
$A$-bimodule $M$ is
     the   differential module
    $$
    {\hC}^\ast  (A;M) =  (\mbox{Hom}   (\cT(sA), M), D_0+D_1) \,.
    $$
  (It is important here to remark
 that the differential graded module   ${\hC}^\ast  (A;M)$   does not coincide with none  of the differential graded modules  $\mbox{Hom} (\overline{\bB}(A), M)$ and 
$\mbox{Hom} (\widetilde{\bB}(A), M)$.)
The {\it Hochschild cohomology of $A$ 
with coefficients in $M$} is
 $$
 \hH^*(A;M)=  H( \hC ^\ast(A;M)) =H((\mbox{Hom}   (\cT(A), M), D_0+D_1))\,.
 $$

 Let $\varphi: A \to A'$  be a homomorphism
 of differential graded algebras.
 Then $A' $ is a differential graded bimodule via $\varphi$ and  $\hC^*(A;A')$ 
is a differential graded algebra.

\vspace{5mm}
\noindent {\bf 3.3}  Let $\psi:C\rightarrow C'$ be a homomorphism of differential 
 graded coalgebras. We denote by $\mbox{Coder}_{\psi}(C,C')$ the 
 subcomplex of $\mbox{Hom}(C,C')$ consisting of 
 $\psi$-coderivations and by 
$
\beta_{C,V} :  \mbox{Hom}(C, V) \rightarrow\mbox{Coder}_{\varphi}(C, \cT (V) )
$ 
the linear isomorphism  extending each linear map  into a coderivation.

 \vspace{3mm}\noindent 
 {\bf Proposition.} (\cite{Ge,St,G:Cartan_Homotopy}) {\sl 

 a) The degree 1 linear isomorphism $\beta_A$:
$$  \hC^*(A;A) =\mbox{Hom}(T\!C(sA), A)  
 \stackrel{\mbox{\small Hom}(T\!C(sA), s)} {\longrightarrow}
 \mbox{Hom}(T\!C'(sA),sA)\stackrel{\beta_{T\!C(sA),sA}}{ \longrightarrow} \mbox{Coder}(\bBnu A)$$
satisfies $ (D_0+D_1) \beta_A = - D \beta_A $,

b) the structure of pre-Gerstenhaber algebra defined
 by
 $\beta_A$ induces the usual structure of Gerstenhaber
algebra on the Hochschild cohomology $HH^{\ast}(A;A)$. }

\vspace{3mm}
Observe that $\beta_A ^{-1}=\mbox{Hom}(T\!C(sA), \tau_{\bBnu A})$
and that for $g\in\mbox{Hom}(\cT(sA),A)$,
$$\beta_A(g)[a_1|\dots|a_n]=
\sum_{0\leq i\leq j\leq n}
(-1)^{\vert s\circ g\vert (\vert sa_1\vert +\dots+\vert sa_i\vert)}
[a_1|\dots|a_i|g([a_{i+1}|\dots|a_j])|a_{j+1}|\dots|a_n].$$

\vspace{5mm}
\noindent {\bf 3.4. Naturality.}  Let $\varphi :  A \to A'$
and $ \psi '  : A' \to B$  be
 homomorphisms of differential graded algebras and $\psi:= \psi' \circ
 \varphi$. Then the two natural maps
  $$
 \hC^\ast(A; A') \stackrel{\hC^\ast (A; \psi')}{\longrightarrow}  \hC^\ast(A; B)
 \stackrel{\hC^\ast ( \varphi ; B)}{\longleftarrow  } \hC^\ast(A'; B)
 $$
 are homomorphisms of differential graded algebras. Let  $\psi'_1$ and
 $\varphi _1$ be the obvious maps  which make commutative the 
 following diagram  of homomorphisms of differential graded modules:
$$
\begin{array}{lclcl}
 \hC^\ast(A; A')& \stackrel{\hC^\ast (A; \psi')}{\longrightarrow }&  \hC^\ast(A;
 B)& \stackrel{\hC^\ast ( \varphi ; B)}{\longleftarrow  }&
 \hC^\ast(A', B)\\
\Phi_{A,A'}\uparrow \cong&&\Phi_{A,B}\uparrow \cong &&\Phi_{A,B}\uparrow \cong \\
\mbox{Hom}_{A^e}(\mathbb B(A;A;A), A')&\stackrel{\psi'_1}{\longrightarrow}&
 \mbox{Hom}_{A^e}(\mathbb B(A;A;A), B)&\stackrel{\varphi_1}{\longleftarrow}&
 \mbox{Hom}_{{A'}^e}(\mathbb B(A';A';A'), B)
\end{array}
$$
It follows from 
(\cite{FHT}-Proposition 2.3) that if $A$, $A'$, $B$ are $\bk$-free
modules and if  $\varphi $ and $\psi$ 
 are quasi-isomorphisms then $\psi_1=\mbox{Hom}_{A^e} ( {\bB}(A;A;A),
 \psi ')$
 and 
$\varphi_1 = \mbox{Hom}_{{A'}^e} ( \varphi , B)$
are quasi-isomorphisms. Thus, so are $ {\bf C}^\ast (A; \psi')$ and ${\bf 
C}^\ast (\varphi; B)$.

 \vspace{3mm}\noindent {\bf Proposition.} {\it   If $f : A \to B$  is a
quasi-isomorphism of differential graded algebras and if $A$ and $B$
are $\bk$-free modules  then the composite,
denoted $ \hH^\ast(f)$
$$
\hH^\ast(A,A) \stackrel{\hH^\ast(A;f)}{\longrightarrow} \hH^\ast(A;B)
\stackrel{(\hH^\ast(f;A))^{-1} }{\longrightarrow} \hH^\ast(B;B)
$$
 is an isomorphism of Gerstenhaber algebras.}

 \vspace{3mm}\noindent {\bf Proof.} As observed above the maps
$$
{\bf C}^\ast(A;A) \stackrel{\small  {\bf C}^\ast (A; f)}\longrightarrow
{\bf C}^\ast(A;B) \stackrel{\small  {\bf C}^\ast (f;A)}\longleftarrow
{\bf C}^\ast(B;B)
$$
are quasi-isomorphisms of differential graded algebras. Let $f_1$ and
$f_2$ be the natural maps which make commutative the following diagram of
homomorphisms of differential graded modules:
$$
\renewcommand{\arraystretch}{1.6}
\begin{array}{lcccl}
 s\hC^\ast (A;A)& \stackrel{ s\hC^\ast(A;f) }{\longrightarrow}& s\hC^\ast
(A;B) &
 \stackrel{s\hC^\ast (f;B) }{\longleftarrow}&  s\hC^\ast (B;B) \\
{\scriptstyle\beta _A} \downarrow \cong && {\scriptstyle \beta _{A,B}} 
\downarrow \cong && {\scriptstyle
\beta _{B}}
\downarrow \cong\\
 \mbox{Coder} (\bBnu A) &\stackrel{f_1}{\longrightarrow}&
\mbox{Coder}_{\bBnu (f)}(\bBnu A,\bBnu B))&
 \stackrel{f_2}{\longleftarrow }& \mbox{Coder} (\bBnu B)
\,,
\end{array}
\renewcommand{\arraystretch}{1}
$$
We have  to prove that the
composite  $H(f_2)^{-1}\circ H( f_1)$
  is a homomorphism of graded Lie algebras. By (\cite{FHT}, Proposition
  3.1) the quasi-isomorphism  $f$
 factors as the composite of
  two quasi-isomorphisms of differential graded algebras
  $$A\build\hookrightarrow_i^{\simeq}
  A\coprod T( V)\build\twoheadrightarrow_p^{\simeq} B$$
  where $i:A\hookrightarrow A\coprod T( V)$ admits a $\bk$-linear retraction
 $r$,  and $p:A\coprod T( V) \twoheadrightarrow B$ admits a $\bk$-linear
  section $s$.
  By functoriality, we have $H\!H^{*}(f )=H\!H^{*}(p)\circ H\!H^{*}(i)$.
  It suffices therefore to prove that $H\!H^*(p)$ and $H\!H^*(i)$ are homomorphisms 
of graded Lie algebras.

In the case $f = p$, $f_1$ admits the $\bk$-linear section
$\mbox{Hom}( T\!C (sA),r)$ and so is surjective.
  Let $x_i$, $i=1,2$, be cycles in $\mbox{Coder}(\bBnu(B))$. Since $f_1$
  is a surjective quasi-isomorphism of complexes,
there exists  cycles, $y_i$, in $\mbox{Coder} \bBnu(A)$ such that
  $$\bBnu(f)\circ y_i=f_1(y_i)=f_2(x_i)=x_i\circ\bBnu(f).$$
  Thus
  $$
\begin{array}{ll}
f_1([y_1,y_2])&=\bBnu(f)\circ y_1\circ y_2
  -(-1)^{\vert y_1\vert\vert y_2\vert}\bBnu(f)\circ y_2\circ y_1\\
&=x_1\circ x_2\circ\bBnu(f)
  -(-1)^{\vert x_1\vert\vert x_2\vert}x_2\circ x_1\circ\bBnu(f)\\
 & =f_2([x_1,x_2])
\,, \end{array}
$$
and if $a_i$ (resp. $b_i$) denotes the class of $x_i$ (resp. $y_i$),
$i=1,2$ then $
sH\!H^*(f)[b_1,b_2] = H(f_2)^{-1}\circ H(f_1) ([b_1,b_2]) =
  [a_1, a_2]\,.$

  In the case $f = i$,
   $f_2$ is surjective and the same argument works,
   mutatis mutandis.

  \hfill  $\square$

\section{Proof of
 Theorem 2.}
In this section $A$ will denote a differential graded algebra
$(TV,d)$ (not necessarly a free model).
The proof of Theorem 2 is a direct consequence of Lemmas 4.2, 4.3 and
4.7 below.

\vspace{3mm}\noindent{\bf 4.1.} Let $\widetilde{\mbox{Der}} A=\mbox{Der} A\oplus sA$ be as defined
in the introduction. We define the injective degree $-1$ linear map
$$i_A:\widetilde{\mbox{Der}} A\hookrightarrow\mbox{Hom}(\cT(sA),A)$$
as follow.
Consider the counity
$\varepsilon_{\cT(sA)}:\cT(sA)\twoheadrightarrow\bk$
and the twisting cochain
$\tau_{\bBnu A}:\cT(sA)\twoheadrightarrow A$.
Given $a\in A$, we denote by $\bar a:\bk\rightarrow A$ the $\bk$-linear map
such that $\bar a(1_{\bk})=a$.
We put
$$\left\{
\begin{array}{ll} 
i_A(sa)=\bar a\circ\varepsilon_{\cT(sA)}&\mbox{ if }a\in A\\
i_A(x)=(-1)^{\vert x\vert}x\circ\tau_{\bBnu A}&\mbox{ if } x\in\mbox{Der}A.
\end{array}
\right.$$

Let $\widetilde{A}=(T(V\oplus \bk\varepsilon),\widetilde{d})$
be the differential graded algebra
where the differential $\widetilde{d}$
is related to $d$ by
$$
 \left\{
 \begin{array}{ll}
 \widetilde{d}\varepsilon = \varepsilon^2, &  \\
 \widetilde{d}v = dv + [\varepsilon ,v] &\mbox{if }v\in V.
 \end{array}
 \right.
 $$
We define 
$$j_{A}  : \widetilde{\rm Der}\, A\to {\rm Der}\,\widetilde{A},
\quad \theta+sa\mapsto j_{A} (\theta + sa) $$
where $j_{A} (\theta + sa)$ is the unique derivation of $\widetilde{A}$
 defined by
 $$
 \left\{
 \begin{array}{l}
j_{A} (\theta + sa)(v) = \theta (v)\,\mbox{for }v\in A\\
j_{A} (\theta+sa)(\varepsilon ) = (-1)^{|a|+1} a.
\end{array}
\right.
$$

\vspace{3mm}\noindent {\bf Lemma 4.2.} 
 {\it The linear maps
$$
s{\hC}^\ast (A;A) \stackrel{s\circ i_A}{\longleftarrow} \widetilde {\rm Der} A 
\stackrel{j_A}\longrightarrow  {\rm Der} \widetilde A
$$
are injective homomorphisms of graded Lie algebras.
}

\vspace{3mm}\noindent{\bf Proof.} 
Precise that the isomorphism
$\beta_{A}\circ s^{-1}:s{\hC}^\ast (A;A)\rightarrow\mbox{Coder}(\bBnu A)$
defines on $s{\hC}^\ast (A;A)$ the Gerstenhaber bracket
$[-,-]$ explicitly defined by:
$$
\left. 
\begin{array}{ll}
(sf)\overline{\circ}(sg) &:=s(f\circ[\beta_A(g)]),\\
 {[} sf, sg{]} &:=(sf)\overline{\circ}(sg)-(-1)^{\vert sf\vert\vert
  sg\vert}(sg)\overline{\circ}(sf)
\end{array}
\right\}\quad f\,, g\in\mbox{Hom}(\cT(sA),A).
$$
From this formula, we deduce now that the inclusion
$s\circ i_A:\widetilde{\mbox{Der}}A\rightarrow s\hC(A;A)$ is a
morphism of graded Lie algebras.
One check directly that $j_A$ is a homomorphism of graded Lie
 algebras.

\hfill$\square$

\vspace{3mm}\noindent {\bf Lemma 4.3.}
{\it  If A is a free model
then 
$$s\circ i_{A}:\widetilde {\rm Der}(A)
\buildrel{\simeq}\over\longrightarrow s\hC^\ast(A;A)$$
is a quasi-isomorphism of differential graded modules.}

\vspace{3mm}\noindent{\bf Proof.}
A straightforward computation shows that the inclusion
$s\circ i_A:\widetilde{\mbox{Der}}A\rightarrow s\hC(A;A)$ is a
morphism of differential graded modules.
In order to prove that $s\circ i_{A}$ is a quasi-isomorphism,
we introduce a quasi-isomorphism between semifree $A^{e}$-modules
$$\Pi:\bB(A;A;A)\buildrel{\simeq}\over\twoheadrightarrow
\widetilde{\mathcal{K}}_{A}$$
and a commutative diagram of differential graded modules
$$
\begin{array}{lll}
\mbox{Hom}_{A^{e}}\left( \widetilde{\mathcal{K}}_{A},A \right) & \build\longrightarrow_{\cong}^{\widetilde{\alpha}_{A}}
&\widetilde {\rm Der}(A)\\
\mbox{Hom}_{A^{e}}\left(\Pi,A \right)\downarrow && \downarrow i_{A}\\
\mbox{Hom}_{A^{e}}\left(\bB(A;A;A),A\right)
&\build\longrightarrow_{\Phi_{A}}^{\cong}
&\hC^\ast(A;A)
\end{array}\quad\quad(4.4)
$$
where $\widetilde{\alpha}_{A}$ is a degree $1$ isomorphism and
$\Phi_{A}$ is the isomorphism defined in 3.2.
Then, by \cite[Proposition 2.3]{FHT},
$\mbox{Hom}(\Pi,A)$ is a quasi-isomorphism and thus so
is $i_A$.
In the remaining of the proof, we precise the definitions of
$\widetilde{\mathcal{K}}_{A}$, $\Pi$ and $\widetilde{\alpha}_{A}$.

Let $S_V$ be the universal derivation
of $A$-bimodules:
$$S_V:A\rightarrow A\otimes V\otimes A\,,
\,\,\, v_1\dots v_n\mapsto\sum_{i=1}^{n}v_1\dots v_{i-1}\otimes
v_i\otimes v_{i+1}\dots v_n\,.$$
Denote by $\overline{S}_V:=(A\otimes s\otimes A)\circ S_V
:A\rightarrow A\otimes sV\otimes A.\quad\quad (4.5)$\\
We put $\widetilde{\mathcal{K}}_{A}:=\left( A \otimes (\bk \oplus sV) \otimes A\,, \widetilde{\delta} \right) $ with 
$\widetilde{\delta} $ defined by $$\left\{
\begin{array}{lll}
\widetilde{\delta} (1\otimes \lambda \otimes 1) & = 0  &\mbox{ if } \lambda \in \bk\\
\widetilde{\delta} (1\otimes sv  \otimes 1) &=
v\otimes 1_{\bk}\otimes 1-\overline{S}_V(dv)-1\otimes 1_{\bk}\otimes v

&\mbox{ if } v  \in V.
\end{array}
\right.
$$

\noindent The map
$\Pi:\bB(A;A;A)\rightarrow\widetilde{\mathcal{K}}_{A}$
is defined by
$\Pi([sa_1\vert\cdots\vert sa_{n}])=\left\{
\begin{array}{ll}
\overline{S}_V(a_1) &\mbox{ if } n=1\\
0 &\mbox{ if } n\neq 1.
\end{array}\right.$
\noindent By~\cite[Theorem 1.4]{Vi} and by~\cite[p. 205 and p. 207]{Luc},
$\Pi$ is
a quasi-isomorphism between semifree $A$-bimodules.
We define the degree $1$ linear isomorphism
$$ \widetilde{\alpha}_{A}:\mbox{Hom}_{A^{e}}\left(\widetilde{\mathcal{K}}_{A},A \right)
\rightarrow \widetilde {\rm Der}(A)={\rm Der}(A)\oplus s(A)\quad\quad (4.6)$$
by
$\left\{
\begin{array}{lll}
\widetilde{\alpha}_{A} (x) &=(-1)^{\vert x\vert} x\circ\overline{S}_V
&\mbox{ if }x\in\mbox{Hom}_{A^{e}}\left( A \otimes sV\otimes
  A,A \right),\\ 
\widetilde{\alpha}_{A} (y) &=-s(y(1\otimes 1_{\bk}\otimes 1))
&\mbox{ if }y\in\mbox{Hom}_{A^{e}}\left( A \otimes \bk\otimes
  A,A \right).
\end{array}\right.$\\
A straightforward computation shows that the diagram (4.4)
commutes.

\hfill $\square$

\vspace{3mm}\noindent {\bf  Lemma 4.7.}
{\it If $A$ is a free model then the injective linear map
$$j_{A}:\widetilde {\rm Der}(A)
\buildrel{\simeq}\over\longrightarrow {\rm Der}(\widetilde{A})$$ is a quasi-isomorphism of differential graded modules.
}

We begin the proof of Lemma 4.7 with a general construction.
Consider the degree $1$ "universal derivation" 
$\overline{S}_V:A=TV\rightarrow A\otimes sV\otimes A$
defined in (4.5).

\vspace{3mm}\noindent {\bf Sub-lemma 4.8.}
{\it

i) $\mathcal{K}_{A}:=(A\otimes sV \otimes A,\delta)$
is a differential $A^{e}$-module with $\delta:A\otimes sV\otimes A\rightarrow A\otimes sV\otimes A,\quad
1\otimes sv  \otimes 1\mapsto -\overline{S}_V(dv),\quad
v\in V.$

ii) The map
$$\alpha_A:\mbox{Hom}_{A^{e}}\left(\mathcal{K}_A,A\right)\rightarrow
{\rm Der}(A),\quad
f\mapsto (-1)^{\vert f\vert} f\circ\overline{S}_V$$
is a degree $1$ isomorphism of differential graded modules.

iii) Denote by $\gamma_{sV}$ the composite
$$\mbox{Hom}(sV,A)\buildrel{\mbox{Hom}(s^{-1},A)}\over\rightarrow
\mbox{Hom}(V,A)\rightarrow {\rm Der}(A)$$
where the map $\mbox{Hom}(V,A)\rightarrow {\rm Der}(A)$ is the
canonical
isomorphism which extends a linear map into a derivation of $A$.
The following diagram of graded modules commutes

\[\xymatrix{
\mbox{Hom}(sV,A)\ar[rr]^{\gamma_{sV}}
&&{\rm Der}(A)\\
&\mbox{Hom}_{A^{e}}(\mathcal{K}_A,A)\ar[ul]^{\Phi_A}
\ar[ur]_{\alpha_{A}}
}\]
where $\Phi_A$ is the canonical isomorphism which restricts an 
$A^{e}$-linear map into a $\bk$-linear map (cf 3.2).
}

\vspace{3mm}\noindent {\bf Proof.}
A straightforward computation proves i) and ii).
In order to prove iii), we use the fact that
$(\gamma_{sV})^{-1}=-\mbox{Hom}(\tau_{A},A)$
where $\tau_{A}:sV\rightarrow A$, $sv\mapsto v$,
$v\in V$.

\hfill$\square$

\vspace{3mm}\noindent {\bf 4.9 Remark.}
An explicit formula for $\gamma_{sV}$ is:
$$\gamma_{sV}(\varphi)=(-1)^{\vert\varphi\vert}
\mu_{A}^{2}\circ (A\otimes (\varphi\circ s)\otimes A)\circ S_V$$
for $\varphi\in\mbox{Hom}(sV,A)$. Here
$\mu_{A}^{2}:A\otimes A\otimes A\rightarrow A$ denotes the
iterated
multiplication of $A$.
The differential $D^{A}$ defined on $\mbox{Hom}(sV,A)$
via the degree $1$ isomorphism $\gamma_{sV}$ is explicitly given by
$$D^{A}(\varphi)=d\circ\varphi+\gamma_{sV}(\varphi)\circ
d_{_{\vert_V}}    \circ s^{-1}
=d\circ\varphi+(-1)^{\vert\varphi\vert}
\mu_{A}^{2}\circ (A\otimes (\varphi\circ s)\otimes A)\circ S_V\circ
d_{_{\vert_V}}    \circ s^{-1}$$
for $\varphi\in\mbox{Hom}(sV,A)$.
In particular, the diagram above is a commutative diagram of
differential graded modules.

\vspace{3mm}\noindent{\bf 4.10. End of proof of Lemma 4.7.}
Let $\widetilde{A}=(T(V\oplus \bk\varepsilon),\widetilde{d})$
defined in 4.1 and consider the following diagram of graded modules
\[
\begin{array}{c}
\xymatrix{
\mbox{Hom}_{\widetilde{A}^{e}}\left(\mathcal{K}_{\widetilde{A}},\widetilde{A} \right) \ar[r]_{\cong}^{\alpha_{\widetilde{A}}}
&{\rm Der}(\widetilde{A})\\
\mbox{Hom}_{A^{e}}\left( \widetilde{\mathcal{K}}_{A},A \right)
\ar[r]_{\cong}^{\widetilde{\alpha}_{A}}\ar@{.>}[u]^{k_A}_{\simeq}
&\widetilde {\rm Der}(A)\ar[u]_{j_A}
}
\end{array}\quad\quad(4.11)
\]
where $\widetilde{\alpha}_{A}$ is the isomorphism of differential
graded modules defined in (4.6)
and $\alpha_{\widetilde{A}}$ is the isomorphism of graded modules
defined in the sublemma 4.8 with $\widetilde{A}=(T(V\oplus
\bk\varepsilon),\widetilde{d})$ in place of $A$.
Our strategy to prove Lemma 4.7 is to complete diagram (4.11)
by a quasi-isomorphism of differential graded modules
$k_A:\mbox{Hom}_{A^{e}}\left( \widetilde{\mathcal{K}}_{A},A \right)
\rightarrow
\mbox{Hom}_{\widetilde{A}^{e}}\left(\mathcal{K}_{\widetilde{A}},\widetilde{A}
\right)$ so that it commutes.

First observe that
$$\mathcal{K}_{\widetilde{A}}=\left(\widetilde{A}\otimes
s(V\oplus \bk\varepsilon)\otimes \widetilde{A},\delta\right)
=\left(\widetilde{A}\otimes(\bk\oplus sV)\otimes
  \widetilde{A},\delta\right).$$
We introduce the differential graded algebra
$$\widehat{A}=(T(V\oplus \bk\varepsilon),\widehat{d})\mbox{ with }
 \left\{
 \begin{array}{l}
 \widehat{d}\varepsilon = -\varepsilon^2 \,,  \\
 \widehat{d}v = dv\,.
 \end{array}
 \right.
 $$
Then the natural inclusion of differential graded algebras 
$$i:A\buildrel{\simeq}\over\hookrightarrow \widehat{A}=A\coprod (T(\varepsilon),\widehat{d})$$ is a
quasi-isomorphism since
$(T(\varepsilon),\widehat{d})$ is acyclic.
Therefore by ~\cite[Proposition 2.3]{FHT},
$\mbox{Hom}_{A^{e}}(\widetilde{\mathcal{K}}_A,i)$ a quasi-isomorphism
of differential graded modules, since $\widetilde{\mathcal{K}}_A$
is a semi-free $A^{e}$-module.
Now we have the following diagram of graded modules
\[
\xymatrix{
\mbox{Hom}_{\widetilde{A}^{e}}\left(\mathcal{K}_{\widetilde{A}},\widetilde{A}
\right)\ar[r]_{\cong}^{\alpha_{\widetilde{A}}}
&{\rm Der}(\widetilde{A})\\
\mbox{Hom}_{\widehat{A}^{e}}\left(\widehat{A}^{e}\otimes_{A ^{e}}
  \widetilde{\mathcal{K}}_{A},\widehat{A}\right)
\ar[u]^{Id}\\
\mbox{Hom}_{A^{e}}\left(
  \widetilde{\mathcal{K}}_{A},\widehat{A}\right)
\ar[u]^{\cong}\\
\mbox{Hom}_{A^{e}}\left( \widetilde{\mathcal{K}}_{A},A \right)
\ar[r]_{\cong}^{\widetilde{\alpha}_{A}}
\ar[u]^{\mbox{Hom}_{A^{e}}\left(\widetilde{\mathcal{K}}_{A},i\right)}_{\simeq}\ar@/^8pc/[uuu]^{k_A}
&\widetilde {\rm Der}(A)\ar[uuu]_{j_A}
}
\]
where the isomorphism of differential graded modules
$$\mbox{Hom}_{A^{e}}\left(\widetilde{\mathcal{K}}_{A},\widehat{A}\right)
\buildrel{\cong}\over\rightarrow
\mbox{Hom}_{\widehat{A}^{e}}\left(\widehat{A}^{e}\otimes_{A ^{e}}
  \widetilde{\mathcal{K}}_{A},\widehat{A}\right)$$
is given by "extension of scalars".
The identity map $Id$ commutes with the differentials. Indeed, the
differential $\widehat{\delta}$ on $\widehat{A}^{e}\otimes_{A ^{e}}
  \widetilde{\mathcal{K}}_{A}$ is given by
$\widehat{\delta}(a)=\delta(a)-[\varepsilon,a]$ for
$a\in\widetilde{A}\otimes(\bk\oplus sV)\otimes
  \widetilde{A}$.
The proof ends by checking that the above diagram commutes.

\hfill $\square$

\section{The Hochschild cochain complex of a differential graded
 coalgebra}

\vspace{5mm}
\noindent{\bf 5.1}  Let $(C,d)$ be a  supplemented
  differential
   graded
   coalgebra, $C = \bk \oplus \overline{C}$, and $(R,d)$  (resp. $(L,d)$)  be a
 right (resp. left)  $C$-comodule. The {\it two-sided
 cobar constructions},
  $\Omega (R; C; L)$  and $\overline{\Omega}(R; C; L)$
    are defined as follows:
    $$\Omega(R;C;L)= (R\otimes T(s^{-1}C)\otimes L\,,d_0+d_1)\,, \quad
\overline{\Omega}(R; C; L)=(R\otimes T(s^{-1}\overline{C})\otimes L,d_0+d_1)\,. 
$$
  A generic element is denoted  $r\langle c_1|c_2|\cdots |c_k\rangle l$ with
 degree $ |r|+|l|+ \sum _{1=1}^k |s^{-1}c_k]$. The
 differential $d_0$ is the unique derivation extending $-sd$  and $d_1$ is given by the formula:
   $$
   \renewcommand{\arraystretch}{1.6}
   \begin{array}{l}
   d_1(r\langle c_1\vert\cdots \vert c_{p-1}\rangle l)=
   -\sum_k(-1)^{\vert r_k'\vert} r'_k\langle x''_k\vert
   c_1\vert\cdots\vert c_{p-1}\rangle l\\
    +
 \displaystyle\sum_{j=1}^{p-1}\sum_i
   (-1)^{\varepsilon_j+\vert c'_{ji}\vert}\, r\langle
   c_1\vert\cdots\vert c'_{ji}\vert c''_{ji}\vert\cdots\vert
   c_{p-1}\rangle l +\sum_j(-1)^{\varepsilon_j}\,
   r\langle c_1\vert\cdots\vert c_{p-1}\vert y'_j \rangle l''_j
   \end{array}
   \renewcommand{\arraystretch}{1}$$
with $\varepsilon_j = \vert r\vert+\vert s^{-1}c_1\vert+\cdots+\vert
   s^{-1}c_{j-1}\vert$, and $\Delta c_j=\sum_i  
c'_{ji}\otimes c''_{ji}$,
   $\Delta r = \sum_k  r'_k\otimes   x''_k $ and  $\Delta l = \sum_j
 y'_j\otimes l''_j $ denote
  the
   non reduced diagonals (resp. reduced diagonals),
   $c_j, c'_{ji}, c'_{ji}, x'_j,y''_i \in C $ (resp. $\overline{C}$),  $l,
 l''_j \in L$ and $ r, r_i' \in R$.

 Hereafter we will use  the {\it normalized} and  the {\it
   non-counital} cobar constructions
 $$
 \overline{\Omega} C=
 \overline{\Omega}(\bk;C;\bk )=\left( T\!C(s^{-1}\overline C), \overline d\right)\,, \mbox{ and }   \widetilde \Omega C=
 \overline{\Omega}(C\oplus \bk )=\left(T\!C(s^{-1} C) ,\tilde  d\right)  \,.
 $$
 and  the twisting cochain 
$ \tau_{\Omeganu C}:C\rightarrow\Omeganu C \,, \quad
 c\mapsto\langle c\rangle $ of $\Omeganu C$. 

\vspace{5mm}
\noindent {\bf Lemma 5.2.} {\it Let $C$ be a locally conilpotent   differential graded coalgebra.  If $C$ is a  free $\bk$-module  then

a) $\overline\Omega C = (T(V) , \overline d) $ is a free model 

b) $\widetilde \Omega C = (T(V\oplus \bk \epsilon ), \tilde d)$ with
 $\tilde  d\epsilon =\epsilon ^2 $ and $ \tilde d v =\overline  dv +\epsilon v - (-1)^{|v|} v \epsilon $, $v \in V$.}

\vspace{2mm}
\noindent{\bf Proof}

a) Write $ C = \bk 1_C \oplus \overline C $, $\, \overline d =\overline d_0+ \overline 
d_1$,  $\,  V=s ^ {-1}\overline C$  and
$$
V(k)=\left\{ 
\begin{array}{ll}
s ^ {-1}\mbox{Ker} \overline{\Delta} ^{(n)} &\mbox{ if  } k= 2n\\
 V(2n) +s^{-1} (\mbox{Ker} \overline{\Delta} ^{(n+1)}\cap \mbox{Ker} \overline d_0 )
&\mbox{ if  } k= 2n+1
\end{array}
\right.
$$
(Here  $\Delta  ^{(0)} = id_C $ and  $V_0 =0$).  Thus $V$ is the union of the 
increasing sequence $\left( V(k)\right)_{k\geq 0}$.
Now observe that if
$c
\in
\mbox{Ker} \overline{\Delta} ^{(n)} $ then  $\overline{\Delta} (c) \in 
\mbox{Ker} \overline{\Delta} ^{(n-1)}\otimes 
\mbox{Ker} \overline{\Delta} ^{(n-1)}$. Thus $d(V(k)) \subset T(V(k-1)$.

b) Write $\epsilon=\langle 1_c\rangle$ and $\widetilde \Omega C = \overline \Omega 
(C\oplus \bk 1_{\bk})= (T(V\oplus\bk\epsilon) ,\tilde d)$ with
$\tilde d= \tilde d_0+\tilde d_1$ satisfies $\tilde d_0 v = \overline
d_0 v $
and $ \tilde d\epsilon = 
\tilde d_1(\langle1_c\rangle)=\langle1_C|1_C\rangle= \epsilon ^2$ and if $v= 
\langle c\rangle$ and $\overline{\Delta}c=\sum_{i} c_i\otimes c'_i$
then $\tilde d_1 v = (-1)^{\vert c\vert}\langle c\vert 1_c\rangle
+\sum _i  (-1)^{\vert c_i\vert}\langle c_i \vert c'_i\rangle + \langle
1_c\vert c \rangle= (-1)^{\vert c\vert}v\epsilon+\overline d_1 v +\epsilon v 
$.
 
\hfill$\square$

\vspace{3mm}
\noindent {\bf 5.3} 
Let $C$ be a differential graded coalgebra and $N$ a differential
graded $C$-bicomodule. For any graded $\bk$-module $V$, we consider
the natural  isomorphism of graded $\bk$-modules
 $$
\Psi _{N,C} :  \mbox{Hom} (N,V))\to \mbox{Hom}_{C^e} ( N ,C\otimes
V\otimes C) \,,
 \hspace{5mm} \Psi_{N,C}(f) = (C\otimes f\otimes C)\circ \Delta_N\,,$$
where $\Delta_N = (C \otimes \Delta_N^r)\circ \Delta_N^l$, $\Delta_N^r:N\rightarrow N\otimes C$ and
$\Delta_N^l:N\rightarrow C\otimes N$ denoting the left and right comodule maps of $N$.
 Here $\mbox{Hom}_{C} ( L , R)$ denotes the graded $\bk$-module  of
 morphisms of $C$-comodules from $ L$ into $ R$ and $C^e =C\otimes C^{op} $ is the
enveloping coalgebra of $C$.
The inverse of $\Psi _{N,C}$ is simply
$\mbox{Hom}(N,\varepsilon_C\otimes V\otimes\varepsilon_C)$.

The isomorphism $\Psi _{N,C} :  \mbox{Hom} (N, T(s^{-1}C))\to \mbox{Hom}_{C^e} ( N , {\Omega}
(C;C;C))$
carries on the graded module $\mbox{Hom}(N,T(s^{-1}C))$ a
differential $D_0+D_1$ explicitely defined by:
$$\begin{array}{ll}
(D_0+D_1)(\varphi)=&
d_{\Omeganu C}\circ\varphi
- (-1)^{\vert\varphi\vert\vert d_N\vert}\varphi\circ d_N\\ &
 -\left(\mu_{\Omeganu C}\circ(\tau_{\Omeganu
     C}\otimes\varphi)\circ\Delta_N^{l}
 -(-1)^{\vert\varphi\vert\vert\tau_{\Omeganu C}\vert}
 \mu_{\Omeganu C}\circ(\varphi\otimes\tau_{\Omeganu
     C})\circ\Delta_N^{r}\right)\end{array}$$
 for $\varphi\in\mbox{Hom}(N,\Omeganu C)$.

  {\it The Hochschild cochain complex of a
differential graded  coalgebra } $C$ with coefficients in the
 $C$-bicomodule $N$  is  the differential graded module~\cite[p. 57]{GS}:
$$
{\bf C}^{\ast} (N;C) =\left( \mbox{Hom} (N, T(s^{-1}C)), D_0+D_1\right)
  \,.
$${\it The Hochschild
 cohomology of the coalgebra $C$ with coefficients in the bicomodule $N$} is
 $$
 \hH^\ast (N;C) = H({\bf C}^\ast (N;C) ) =H\left(  \mbox{Hom} (N, T(s^{-1}C)),
 D_0+D_1\right)\,.
 $$

 If $\psi : C' \to C$ is a homomorphism of differential
graded coalgebras, then ${\bf C}^\ast (C';C)$ is a
 differential graded algebra.

 \vspace{3mm}\noindent  {\bf 5.4
  Proposition 5.4.} {\sl The degree 1  linear isomorphism  $\gamma_C$:
$$
\hC^*(C;C) = \mbox{Hom} (C, T(s^{-1}C)) \stackrel{\mbox{\small Hom}
(s,   T(s^{-1}C))}\longrightarrow 
         \mbox{Hom} (s^{-1}C, T(s^{-1}C){\rightarrow} \mbox{Der}( \widetilde{\Omega} C)$$
satisfies $\gamma_C\circ (D_0+D_1)=-D\circ\gamma_C$.}

\vspace{3mm}
\noindent{\bf Proof.}
By sub-lemma 4.8 and Remark 4.9 applied to $\Omeganu C$,
the degree $1$ linear map $\gamma_C$:
$$
\left(\mbox{Hom}(C, Ts^{-1}C), D^{\Omeganu C}\right)
\rightarrow \mbox{Der}(\Omeganu C)
$$
anticommutes with the differentials.
A straightforward computation shows that the differential
$D^{\Omeganu C}$ coincides with differential $D_0+D_1$ of
$\hC^*(C;C)$.
 
\hfill$\square$


\noindent{\bf 5.5}  We now suppose that $C$ is a locally conilpotent.
  {\it The normalized Hochschild cochain complex of} $C$ is  the
  differential graded module:
$$
\overline{\hC}^{\ast} (C;C) =\left( \mbox{Hom} (C, T(s^{-1}\overline{C})), D_0+D_1\right)
$$
(The isomorphism $
\Psi _{C} :  \overline{\hC}^{\ast} (C;C)\buildrel{\cong}\over\rightarrow\mbox{Hom}_{C^e} ( C , \overline{\Omega}
(C;C;C))$ defined above 
being an homomorphism of differential graded modules).

\vspace{3mm}
 \noindent{\bf Proposition 5.6.} 
{\it The inclusion $ i : \overline\Omega C \hookrightarrow \widetilde
  \Omega C$ induces a  quasi-isomorphism  of differential  graded algebras 
$$\mbox{Hom}(C,i):
\overline{\hC}^{\ast}(C;C)
\buildrel{\simeq}\over\hookrightarrow \hC^{\ast}(C;C)
$$
 Moreover,
there exists   a  structure of differential graded Lie algebra on
$s\overline{\hC}^{\ast}(C;C)$ such that
$$s\mbox{Hom}(C,i):
s\overline{\hC}^{\ast}(C;C)
\buildrel{\simeq}\over\hookrightarrow s\hC^{\ast}(C;C)$$
is a homomorphism  of differential graded Lie algebras.}

\vspace{3mm}
\noindent{\bf Proof.}
By diagram (4.11) applied when $TV=\overline{\Omega}C$, we obtain 
the following diagram of differential graded modules
\[
\begin{array}{lll}
s\mbox{Hom}_{\Omeganu C^{e}}(\mathcal{K}_{\Omeganu C},\Omeganu C)
&
\build\longrightarrow_{\cong}^{\mbox{Hom}(\overline{S}_{s^{-1}C},\Omeganu
  C)\circ s^{-1}}
&\mbox{Der}(\Omeganu C)\\
sk_{\overline{\Omega} C}
\uparrow&&\uparrow j_{\overline{\Omega}  C}\\
s\mbox{Hom}_{\overline{\Omega}C^{e}}(\widetilde{\mathcal{K}}_{\overline{\Omega}C},\overline{\Omega}C)&
\build\longrightarrow_{\cong}^{\widetilde{\alpha}_{\overline{\Omega}  C}\circ s^{-1}}
&\widetilde{\mbox{Der}}(\overline{\Omega} C)
\end{array}\,.
\]
On the other hand, the following diagram of graded modules commutes obviously
\[
\begin{array}{lllll}
\mbox{Hom}_{C^{e}}(C,\Omega(C;C;C))
&\build\longleftarrow_\cong^{\Psi_C}
&\hC^{\ast}(C;C)
&\build\longleftarrow_\cong^{\Phi_{\Omeganu C}}
&\mbox{Hom}_{\Omeganu C^{e}}\left(\mathcal{K}_{\Omeganu C},\Omeganu
  C\right)\\
\mbox{Hom}_{C^{e}}(C,C\otimes i\otimes C)\uparrow
&&\mbox{Hom}(C,i)\uparrow
&&
k_{\overline{\Omega}C}\uparrow\\
\mbox{Hom}_{C^{e}}(C,\overline{\Omega}(C;C;C))
&\build\longleftarrow_\cong^{\Psi_C}
&\overline{\hC}^{\ast}(C;C)
&\build\longleftarrow_\cong^{\Phi_{\overline{\Omega} C}}
&\mbox{Hom}_{\overline{\Omega}C^{e}}\left(\widetilde{\mathcal{K}}_{\overline{\Omega} C},\overline{\Omega}C\right)
\end{array}\,.
\]
Since the maps  $\Psi_C$'s, $\Phi_{\Omeganu C}$,
$\mbox{Hom}_{C^{e}}(C,C\otimes i\otimes C)$
and $k_{\overline{\Omega}C}$
are homomorphisms of differential graded modules, then
$\Phi_{\overline{\Omega} C}$ and 
$\mbox{Hom}(C,i)$
are also homomorphisms of differential graded modules.
Define $\overline{\gamma}_C:=
\widetilde{\alpha}_{\overline{\Omega}C}\circ\Phi_{\overline{\Omega}C}^{-1}$.
The remaining
of the statement follows also from the commutativity of above
diagrams.

\hfill$\square$

   \section{Proof of Propositions B and C}

\noindent   {\bf  6.1 Proof of Proposition C.}
Let $C$ be a supplemented
conilpotent differential graded coalgebra.
We denote by
$\sigma_C: C\stackrel  {\simeq}\rightarrow    \overline{\bB}\,\overline{\Omega}
  C\stackrel  {\simeq}\hookrightarrow \bB(\bk;\overline{\Omega} C;\bk) $
the counity of the bar-cobar adjunction, $$\sigma_C(c) =[<c>] +
\sum_{i\geq 1}\sum_j [<c_{1,j}> \vert\cdots\vert
  <c_{i+1, j}>]\,,$$
  where $\overline\Delta^ic = \sum_j c_{1,j}\otimes \ldots \otimes
  c_{i+1,j}$, $c\in \overline C$.
Define ${\mathcal D}_2 = \mbox{Hom}(\sigma_C ;\overline{\Omega} C)  : {\bf C}^\ast (\overline{\Omega} C ;\overline{\Omega} C)
{\to }\overline{\bf C}^\ast ( C ;C )$.

\noindent We will  prove that{\it 

 a) $\mathcal {D}_2 $ is a quasi-isomorphism of differential graded
 algebras
which admits a section $\Gamma$.

 b) The map $ s\Gamma :  s  \overline{\bf C}^\ast ( C  ;C
)\to s {\bf C}^\ast ( \overline{\Omega} C ;\overline{\Omega} C ) $
  is a homomorphism of differential graded Lie algebras.}

\vspace{3mm}

 a) The quasi-isomorphism $
\Pi :\bB(\overline{\Omega}C;\overline{\Omega}C;\overline{\Omega}C)
\rightarrow\widetilde{\mathcal{K}}_{
\overline{\Omega}C}$ of semifree
 $\overline{\Omega}C^{e}$-modules
described in the proof of Lemma 4.3  with $A=\overline{\Omega}C$,
admits the section
$\nabla_\infty:=\overline{\Omega}C\otimes\sigma_C\otimes\overline{\Omega}C$
which is a homomorphism of differential graded
modules~\cite[p. 209]{Luc} and we have the following commutative
diagram
of graded modules
 $$
 \begin{array}{ccc}
\left(\mbox{Hom}(C,\overline{\Omega }C),D_0+D_1\right)
&\build\longleftarrow_{\cong}^{\Phi_{\overline{\Omega} C}}
&\mbox{Hom}_{\overline{\Omega}C^{e}}
(\widetilde{\mathcal{K}}_{\overline{\Omega}C},\overline{\Omega }C)\\
\hspace{-4mm}\mbox{Hom}(\sigma_C,\overline{\Omega}C)\uparrow
&&
\uparrow\mbox{Hom}_{\overline{\Omega}C^{e}}
(\nabla_\infty,\overline{\Omega}C)\\
\left(\mbox{Hom}(\bBnu\overline{\Omega }C,\overline{\Omega }C),D_0+D_1\right)
&\build\longleftarrow_{\cong}^{\Phi_{\overline{\Omega} C}}
&\mbox{Hom}_{\overline{\Omega}C^{e}}
(\bB(\overline{\Omega}C;\overline{\Omega }C;\overline{\Omega }C),
\overline{\Omega }C)
 \end{array}
 $$
where the upper homomorphism is defined as the lower
homomorphism has been defined in 3.2.
Since the right vertical  map is a quasi-isomorphism of differential
graded modules and since the horizontal
maps are isomorphisms of differential graded modules, the homomorphism
$\mathcal{D}_2:=\mbox{Hom}(\sigma_C,\overline{\Omega}C)$ is also a quasi-isomorphism
of differential graded modules.
The section $\Gamma$ is then defined by
$$\Gamma:=\Phi_{\overline{\Omega} C}\circ
\mbox{Hom}_{\overline{\Omega}C^{e}}\left(\Pi ,\overline{\Omega}C
  \right)\circ\Phi_{\overline{\Omega} C}^{-1}.
$$

b) By Lemma 4.2, $$s\circ i_{\overline{\Omega}C}:
\widetilde{\mbox{Der}}\overline{\Omega}C\hookrightarrow
s {\bf C}^\ast ( \overline{\Omega} C ;\overline{\Omega} C )$$ is a
homomorphism of differential graded Lie algebras and diagram (4.4)
becomes 
\[
\begin{array}{ccc}
s\mbox{Hom}_{\overline{\Omega}C^{e}}
\left(\widetilde{\mathcal{K}}_{\overline{\Omega}C},
\overline{\Omega}C \right) & \build\longrightarrow_{\cong}^{\widetilde{\alpha}_{\overline{\Omega}C}\circ s^{-1}}
&\widetilde {\rm Der}(\overline{\Omega}C)\\
\hspace{-9mm} s\mbox{Hom}_{\overline{\Omega}C^{e}}\left(\Pi ,\overline{\Omega}C
  \right)\downarrow && \downarrow s\circ i_{\overline{\Omega}C}\\
s\mbox{Hom}_{\overline{\Omega}C^{e}}\left(\bB(\overline{\Omega}C;\overline{\Omega}C;\overline{\Omega}C),\overline{\Omega}C\right)
&\build\longrightarrow_{s\Phi_{\overline{\Omega}C}}^{\cong}
&s {\bf C}^\ast ( \overline{\Omega} C ;\overline{\Omega} C ) 
\end{array}\]
Therefore $s\Gamma$ is a homomorphism of differential graded Lie algebras.

\hfill $\square$

\vspace{3mm}
 \noindent {\bf  6.2 Proof of Proposition B.}

  Let $C$ be a differential graded  coalgebra   and consider the
 linear map $$
 \mathcal {D}_1 :  {\bf C}^\ast (C;C) \to  {\bf C}^\ast (C^\vee
 ;C^\vee )\,, \quad  f \mapsto  f^\vee \circ \Theta
  $$
where  $\Theta :\cT(sA)\rightarrow(T(s^{-1}C))^{\vee}$ is defined by:
 $$
\Theta([f_1\vert\cdots\vert f_n])(\langle c_1\vert\cdots\vert c_k\rangle )
 =
\left\{ \begin{array}{ll}
 0 &  \mbox{if }k\neq n,\\
 (-1)^{n}\varepsilon_\sigma
 f_1(c_1)\cdots f_n(c_n) &\mbox {if }k= n\,.
\end{array}
 \right.
$$
 Here $\varepsilon_\sigma$ is the graded signature obtained by the
 strict application of the Koszul rule to the graded permutation\\
 $
s, f_1, s, f_2,\dots, s, f_n, s^{-1}, c_1, s^{-1}, c_2,\dots,
s^{-1},c_n \mapsto 
f_1, s, s^{-1}, c_1, f_2, s, s^{-1}, c_2,\dots, f_n, s, s^{-1}, c_n
$.

The proof decomposes in three steps.

\vspace{3mm}\noindent {\bf 6.3 Step 1.}  {\sl  The  map
  $\mathcal {D}_1$ is a   homomorphism of   differential graded modules.
Moreover, 

\noindent i) $\mathcal {D}_1$ is an isomorphism whenever
 $C=C^{\geq 0}$ is a free graded $\bk$-module of finite type.

\noindent ii)  $\mathcal {D}_1$ is a quasi-isomorphism whenever $C$ is coaugmented and
 $\overline{C}=\overline{C}_{\geq 2}$ is a
free graded $\bk$-module of finite type.}

Let $A=C^\vee$ and consider  the following commutative diagram of graded
$\bk$-modules.
$$\xymatrix{
\mbox{Hom}_{C^{e}}(N,\Omega(C;C;C))
\ar[r]^{t}\ar[d]_{\Psi_{N,C}^{-1}}^{\cong}
&\mbox{Hom}_{A^{e}}(\Omega(C;C;C)^{\vee},N^{\vee})
\ar[r]^{{\scriptstyle \mbox{Hom}(\Theta_C,N^{\vee})}}\ar[d]
&\mbox{Hom}_{A^{e}}(\bB(A;A;A),N^{\vee})\ar[d]^{\Phi_{A,N^{\vee}}}_{\cong}\\
\left(\mbox{Hom}(N,T(s^{-1} C)),D_0+D_1 \right)\ar[r]_{\overline{t}}\ar@{=}[d]
&\mbox{Hom}((T (s^{-1}C))^{\vee},N^{\vee})\ar[r]_{{\scriptstyle \mbox{Hom}(\Theta,N^{\vee})}}
&\left(\mbox{Hom}(T\!C(sA),N^{\vee}),D_0+D_1 \right)\ar@{=}[d]\\
{\bf C}^\ast (N;C)&&{\bf C}^\ast (A;N^\vee)}
$$
where $t$ and $\overline{t}$ are the transposition
maps $f\mapsto f^\vee  $ and $\Theta_C:A\otimes T
(sA)\otimes A\rightarrow\left(C\otimes
  T(s^{-1}C)\otimes C\right)^{\vee}$ is the unique homomorphism of
  $A$-bimodules
extending $\Theta$. A tedious computation proves that the maps
$\Theta:\bBnu A\rightarrow(\Omeganu C)^{\vee}$ and $\Theta_C:
\bB(A;A;A)\rightarrow (\Omega (C;C;C))^{\vee} $ commute with the
 differentials. Since the maps on the upper line of the diagram are
 homomorphism of differential graded modules  so is the composite
$\mathcal{D}_1$  of
 the maps appearing in the lower line. 

i) Suppose now that $C$ is free of finite type.  Since $ C= C^{\geq 0}$,
 $s^{-1}C=(s^{-1}C)^{\geq 1}$,
 $T(s^{-1}C)$ is also free of
finite type. Therefore  $t$ is an isomorphism.
Since $sA\cong (s^{-1}C)^{\vee}$, $\Theta$ is also an isomorphism.

ii) Let $A$ be an augmented differential graded algebra such that
$\overline{A}$ is $\bk$-free.
Define the complex
$$\overline{\hC}^{\ast}(A;A):=(\mbox{Hom}(\cT(s\overline{A},A)),D_0+D_1)$$
such that the canonical isomorphism
$\overline{\Phi}_A:\mbox{Hom}_{A^{e}}(\overline{\bB}(A;A;A),A)
\rightarrow \overline{\hC}^{\ast}(A;A)$
is a morphism of differential graded modules.
Let $p_A:\bB(A;A;A)\twoheadrightarrow \overline{\bB}(A;A;A)$
and $p:\bBnu A\twoheadrightarrow\overline{\bB}A$ be the canonical
projections.
Since $p_A$ is a quasi-isomorphism of $A^{e}$-semifree modules and
since the
following
diagram of differential graded modules commutes,
\[
\begin{array}{lll}
\mbox{Hom}_{A^{e}}(\overline{\bB}(A;A;A),A)
&\build\rightarrow_{\simeq}^{Hom_{A^{e}}(p_A,A)}
&\mbox{Hom}_{A^{e}}(\bB(A;A;A),A)\\
\overline{\Phi}_A\downarrow
&
&\downarrow \Phi_A\\
\overline{\hC}^{\ast}(A;A)
&\buildrel{Hom(p,A)}\over\rightarrow
&\hC^{\ast}(A;A)
\end{array}
\]
then $Hom(p,A)$ is a quasi-isomorphism.

Suppose now that $C$ is supplemented and
$\overline{C}=\overline{C}_{\geq 2}$ is $\bk$-free of finite type.
Since $A:={C}^{\vee}$ is $\bk$-free, then $Hom(p,C^{\vee})$ is a
quasi-isomorphism.
Let $i:\overline{\Omega}C\hookrightarrow\Omeganu C$ be the canonical
inclusion.
Consider the unique morphism of graded modules
$\overline{\Theta}:\overline{\bB}C^{\vee}
\rightarrow(\overline{\Omega}C)^{\vee}$
such that $\overline{\Theta}\circ p=i^{\vee}\circ\Theta$
and the linear map
$$
 \overline{\mathcal {D}}_1 :  \overline{\bf C}^\ast (C;C) \to
\overline{\bf C}^\ast (C^\vee
 ;C^\vee )\,, \quad  f \mapsto  f^\vee \circ \overline{\Theta}.
  $$
We have the commutative diagram of graded modules
\[
\xymatrix{
\overline{\bf C}^\ast (C;C)\ar[r]^{\overline{\mathcal {D}}_1}_{\cong}
\ar@{^{(}->}[d]_{\mbox{Hom}(C,i)}^{\simeq}
&\overline{\bf C}^\ast (C^\vee;C^\vee )
\ar@{^{(}->}[d]^{\mbox{Hom}(p,C^{\vee})}_{\simeq}\\
{\bf C}^\ast (C;C)\ar[r]^{\mathcal {D}_1}_{\simeq}
&{\bf C}^\ast (C^\vee;C^\vee )
}
\]
Since the vertical maps are injective, $\overline{\mathcal {D}}_1$
is a morphism of differential graded modules.
As in case i), $\overline{\mathcal {D}}_1$ is an isomorphism and
so using Proposition A, $\mathcal {D}_1$ is a quasi-isomorphism.

 \hfill $\square$

\vspace{3mm}\noindent {\bf 6.4 Step 2.} {\it The linear map $
 \mathcal {D}_1 :  {\bf C}^\ast (C;C) \to {\bf C}^\ast (C^\vee ;{C}^\vee )
 $ is a homomorphism of  graded algebras.}

  Observe that,  without finite type hypothesis,
 $ \Theta: T\!C(sA)\rightarrow \left(T(s^{-1}C)\right )^{\vee}$ is not a homomorphism of
 graded coalgebras. Nevertheless, the composite,
$$\widehat \Theta:
T(s^{-1} C)   \hookrightarrow \left( T(s^{-1} C) \right)^{\vee\vee }
\stackrel{ \left( \Theta \right) ^{\vee } } \rightarrow (T\!C(s
A))^{\vee}
$$
 is  a homomorphism of graded algebras. From the commutative diagram
$$
\renewcommand{\arraystretch}{1.4}
\begin{array}{llllll}
\mbox{Hom}(C,T(s^{-1}C))&\stackrel {\overline t}{\to} &
 \mbox{Hom}((T (s^{-1} C))^{\vee},{C}^{\vee}) \\
\hspace{-6mm} {\scriptstyle \mbox{Hom}({C'},\widehat \Theta )}\downarrow && \hspace{10mm}
\downarrow {\scriptstyle \mbox{Hom}(\Theta,{C}^{\vee})} \\
\mbox{Hom}(C,(T\!C(s A))^{\vee})& \stackrel{\tau}{\to} &
\mbox{Hom}(T\!C(s A), C^{\vee})\,,
\end{array}
\renewcommand{\arraystretch}{1}
$$
where $\tau ( f)$ is the composite $ T\!C(s A ) \hookrightarrow
(T\!C(s A ) )^{\vee\vee} \stackrel{f^\vee}{\to} C^{\vee}$, we
deduce that $\mathcal{D}_1:= \mbox{Hom}(\Theta ,{C}^\vee)\circ
\overline t=\tau\circ \mbox{Hom} (C,\widehat\Theta)$ is a
morphism of graded algebras.

 \hfill $\square$

\vspace{3mm}\noindent {\bf 6.5 Step 3.} {\sl  The homomorphism
 $
 s\mathcal {D}_1 : s {\bf C}^\ast (C;C) \to  s{\bf C}^\ast
(C^\vee
 ;{C}^\vee )
 $
 is an homomorphism of   graded  Lie algebras.}

 We want  to show that
$$\beta_A\circ \mathcal{D}_{1}\circ\gamma_C^{-1}:
 \mbox{Der}(\Omeganu  C)\stackrel{{\gamma_C}^{-1}}{\rightarrow}
{\bf C}^\ast (C;C) \stackrel { \mathcal{D}_1}{\to}
 {\bf C}^\ast (A;A) \stackrel{\beta_A}{\to} \mbox{Coder}(\bBnu A)$$
satisfies 
 $$
(\beta_A\circ \mathcal{D}_{1}\circ\gamma_C^{-1})([\gamma_C(f),\gamma_C(g)])=
 [(\beta_A\circ \mathcal{D}_{1})(f),(\beta_A\circ \mathcal{D}_{1})(g)]\,,
\quad f,g\in\mbox{Hom}(C,\Omeganu C)\,.\quad(6.6)
$$

We  consider the twisting cochains $
\tau_{\Omeganu  C} : C\rightarrow  \Omeganu ( C) $ and $ 
\tau_{\bBnu A}: \bBnu A\rightarrow A  $. By definition, $\gamma_C(f)$   is the
unique derivation such that
 $f = (-1)^{\vert f\vert}
\gamma_C(f)\circ\tau_{\Omeganu  C}$,  and  $\beta_A(g)$  is the
unique  coderivation  such that
 $\tau_{\bBnu A}\circ\beta_A(g)=g$. On the other hand, we have
 clearly
 $\tau_{\bBnu A}=(\tau_{\Omeganu  C})^{\vee}\circ\Theta$.
Momentarily assuming that
the next diagram is commutative
(See Lemma 6.7 below).
$$
\xymatrix{
\bBnu A\ar[r]^{\Theta}\ar@{.>}[d]|-{\beta_A(g^{\vee}\circ\Theta)}
& (\Omeganu C)^{\vee}\ar[d]^{-\gamma_C(g)^{\vee}}\\
\bBnu A\ar[r]^{\Theta}\ar@{.>}[d]|-{\beta_A(f^{\vee}\circ\Theta)}
& (\Omeganu C)^{\vee}\ar[d]^{-\gamma_C(f)^{\vee}}\ar@/^/[ddr]^{f^{\vee}}\\
\bBnu A\ar[r]^{\Theta}\ar@/_/[drr]_{\tau_{\bBnu A}}
& (\Omeganu C)^{\vee}\ar[dr]|-{(\tau_{\Omeganu C})^{\vee}}\\
& & A
}
$$
We have $[\gamma_C(f),\gamma_C(g)]^{\vee}=
(-1)^{\vert  \gamma_C(f) \vert\vert \gamma_C(g)  \vert}[\gamma_C(g)^{\vee},\gamma_C(f)^{\vee}]=-[\gamma_C(f)^{\vee},\gamma_C(g)^{\vee}]
$.
Therefore
the following diagram commutes
$$
\xymatrix{ \bBnu
A\ar[rr]^{\Theta}\ar[d]_{[\beta_A(f^{\vee}\circ\Theta),
\beta_A(g^{\vee}\circ\Theta)]} && (\Omeganu
C)^{\vee}\ar[d]^{-[\gamma_C(f),\gamma_C(g)]^{\vee}}
\ar@/_/[ddrr]^{(\gamma_C^{-1}[\gamma_C(f),\gamma_C(g)])^{\vee}}\\
\bBnu A\ar[rr]^{\Theta}
\ar@/_/[drrrr]_{\tau_{\bBnu A}}
&& (\Omeganu C)^{\vee}\ar[drr]_{\tau_{\Omeganu C}^{\vee}}\\
&& && A.
}
$$
This
implies that the
coderivation
$[\beta_A(f^{\vee}\circ\Theta),\beta_A(g^{\vee}\circ\Theta)]$
coincides with the coderivation $
\beta_A\left((\gamma_C^{-1}[\gamma_C(f),\gamma_C(g)])^{\vee}\circ\Theta\right)
$. Therefore (6.6) is proved since ${\mathcal D}_1(h) = h^\vee
\circ \Theta$.

\hfill $\square$

\vspace{3mm} \noindent{\bf Lemma 6.7.} {\sl Let $f : C \to
\Omeganu C$ be a linear map and $A = C^\vee$. With the
previous notations,
then the next diagram commutes
$$
\xymatrix{
\bBnu A\ar[r]^{\Theta}\ar@{.>}[d]|-{\beta_A(f^{\vee}\circ\Theta)}
& (\Omeganu C)^{\vee}\ar[d]^{-\gamma_C(f)^{\vee}}\\
\bBnu A\ar[r]^{\Theta}
& (\Omeganu C)^{\vee}
}.
$$}

\vspace{3mm}\noindent {\bf Proof.} Recall
the extension maps $\beta_A : \mbox{Hom} (\bBnu A, A) \to
\mbox{Coder} (\bBnu A)$  (3.5), $\gamma_C : \mbox{Hom}
(C , \Omeganu C)\to\mbox{Der}(\Omeganu C) $ (4.2) and if $A=TV$  the universal derivation $S_V$ (4.3). Let $S^V$ be the universal coderivation   of
$\cT (V)$-bicomodules:
 $S^{V}:\cT (V)\otimes V\otimes\cT (V)\rightarrow\cT (V)\,,
\,\,\,\alpha\otimes v\otimes\beta\mapsto \alpha v\beta\,.$ The
morphisms $S^{(V^\vee)}$ and $(S_V)^\vee$ are related by the
following diagrams where $i_1$ and $i_2$ denote standard maps.
$$ \renewcommand{\arraystretch}{1.5}
\begin{array}{ccc}
\cT(V^\vee) \otimes V^\vee\otimes \cT(V^\vee) &
\stackrel{i_1}{\rightarrow} & \left(TV\otimes V\otimes TV\right)^\vee\\
{\scriptstyle S^{(V^\vee)}}\downarrow && \downarrow {\scriptstyle
(S_V)^\vee}\\
\cT(V^\vee) &\stackrel{i_2}{\rightarrow }& (TV)^\vee \end{array}
\renewcommand{\arraystretch}{1}
$$
Then $\gamma_C(f)$ (resp. $\beta_A(g)$) is defined as the
composite
$$\gamma_C(f) : \Omeganu C \stackrel{S_{s^{-1}C}}{\to}
\Omeganu C\otimes s^{-1}C\otimes \Omeganu C \stackrel{1\otimes
((-1)^{\vert f\vert}f\circ s)\otimes 1}{\to} (\Omeganu C)^{\otimes^3}\stackrel{\mu^2}{\to}
\Omeganu C$$ $$(\mbox{resp.}\quad\beta_A(g) : \bBnu A \stackrel{\Delta^2}{\to}
(\bBnu A)^{\otimes^3} \stackrel{1\otimes (s\circ g)\otimes 1}{\to} \bBnu A
\otimes sA\otimes \bBnu A \stackrel{S^{sA}}{\to} \bBnu A).$$ The
claim follows then from the commutativity of the following
diagram
$$ \renewcommand{\arraystretch}{1.5}
\begin{array}{rcl}
\bBnu A & \stackrel{\Theta}{\to} & (\Omeganu C)^\vee\\
{\scriptstyle \Delta^2_{\bBnu A}}\downarrow && \downarrow
{(\mu^2)^\vee}\\
(\bBnu A)^{\otimes^3} & \stackrel{\Theta^{\otimes^3}}{\to}& \left(
(\Omeganu C)^{\otimes^3}\right)^\vee\\
{\scriptstyle 1\otimes( f^\vee \circ \Theta) \otimes 1}\downarrow
&& \downarrow {\scriptstyle (1\otimes f\otimes 1)^\vee}\\
\bBnu A\otimes A\otimes \bBnu A &\stackrel{\Theta \otimes 1
\otimes \Theta}{\to} &(\Omeganu C\otimes C \otimes \Omeganu
C)^\vee\\
{\scriptstyle S^{sA}\circ (1\otimes s\otimes 1)}\downarrow  && \downarrow{\scriptstyle
-(S_{s^{-1}C})^\vee\circ (1\otimes s\otimes 1)^\vee}\\
\bBnu A &\stackrel{\Theta}{\to} & (\Omeganu C)^\vee
\end{array}
\renewcommand{\arraystretch}{1}$$

\hfill $\square$

\vspace{1cm}

\hspace{-1cm}\begin{minipage}{19cm}
\small
\begin{tabular}{lll}
felix@math.ucl.ac.be              &luc.menichi@univ-angers.fr
&jean-claude.thomas@univ-angers.fr\\
D\'epartement de math\'ematique  &  D\'epartement de math\'ematique  &
D\'epartement de math\'ematique \\

 Universit\'e Catholique de Louvain  &Facult\'e des Sciences  &  Facult\'e des
Sciences \\

 2, Chemin du Cyclotron           &2, Boulevard Lavoisier & 2, Boulevard
Lavoisier     \\

 1348 Louvain-La-Neuve, Belgium       & 49045 Angers, France &
 49045
Angers, France

\end{tabular}

\end{minipage}

 \end{document}